\documentclass{amsart}

\usepackage[english]{babel}
\usepackage[latin1]{inputenc}
\usepackage{amsmath}
\usepackage{amsthm}
\usepackage{amssymb}
\usepackage{tikz}
\usepackage{imakeidx}
\usepackage{appendix}
\usepackage{tikz-cd}
\usepackage{verbatim}
\usepackage{enumitem}
\usepackage{multicol}
\usepackage{hyperref}
\usepackage{cleveref}
\usepackage{mathtools}

\newtheorem{thm}{Theorem}[section]

\newtheorem{prop}[thm]{Proposition}
\newtheorem{lem}[thm]{Lemma}

\newtheorem{cor}[thm]{Corollary}

\newtheorem{conj}[thm]{Conjecture}

\theoremstyle{definition}
\newtheorem{defin}[thm]{Definition}

\theoremstyle{remark}

\numberwithin{equation}{section}

\newcommand{\Q}{\mathbb Q}
\newcommand{\Z}{\mathbb Z}
\newcommand{\G}{\mathbb G}
\renewcommand{\L}{\mathbb L}

\newcommand{\Spec}{\operatorname{Spec}}
\renewcommand{\c}{\subseteq}
\newcommand{\A}{\mathbb A}
\newcommand{\mc}[1]{\mathcal{#1}}
\newcommand{\cl}{\overline}
\newcommand{\set}[1]{\{#1\}}
\renewcommand{\phi}{\varphi}
\newcommand{\on}[1]{\operatorname{#1}}

\title{Motivic classes of classifying stacks of some semi-direct products}


\author{Ivan Martino}
\address{Department of Mathematics\\
        Royal Institute of Technology\\
        Stockholm\\ Sweden.}
\email{imartino@kth.se}
\thanks{Ivan Martino has been partially
        supported by the Swiss National Science Foundation Professorship
grant PP00P2\_150552/1; he has been partially supported by the Zelevinsky
Research Instructor Fund.
        Currently, he is supported by the Knut and Alice Wallenberg
Fundation and by the Royal Swedish Academy of Science.}

\author{Federico Scavia}
\address{Department of Mathematics\\
        University of British Columbia\\
        Vancouver, BC \\Canada.} \email{scavia@math.ubc.ca}
\thanks{Federico Scavia was partially supported by a graduate fellowship
from the University of British Columbia.}

\begin{document}
        \begin{abstract}
                \noindent
                Let $k$ be a field, let $G$ be a finite group, and let $T$
be a split $k$-torus on which $G$ acts multiplicatively. For every $m\geq
1$ denote by $T[m]$ the $m$-torsion subgroup of $T$. Under a suitable
assumption on $m$, we show that the motivic class of $B(T[m]\rtimes G)$ in
$K_0(\on{Stacks}_k)$ equals that of $BG$. As a consequence, we prove that
the motivic class of $BW$ is trivial for a large class of complex
reflection groups $W$.
        \end{abstract}

        \maketitle

        \section{Introduction}

        Let $k$ be a field, and let $G$ be a finite group scheme over $k$.
This paper contributes to the computation of the motivic class $\{BG\}$ of
the classifying stack of $G$ in the Grothendieck ring of algebraic stacks
$K_0(\on{Stacks}_k)$. More specifically, we are interested in the
triviality of $\{BG\}$, following up on work of Ekedahl
\cite{ekedahl2009geometric, ekedahl2009grothendieck}, the computations of
the first author \cite{Martino-TEIFFG}, and the results of
\cite{bergh2015motivic, MR3466822, MR3742448, totaro2016motive}.

        Before we go any further, we provide the context where our results
are relevant and several motivations to tackle this motivic computations.

        \subsection*{The Noether Problem}
        Let $V$ be a faithful $G$-representation, finite-dimensional over
$k$. The Noether Problem for $G$ and $V$ is the question of whether the
quotient variety $V/G$ is $k$-rational, that is, birational to some affine
space over $k$. In 1917, Noether \cite{Noether1917} studied this question
in the case when $G$ is constant (i.e. a finite group). In this case, the
problem amounts to determine whether the field of invariants $k(V)^G$ is
rational (i.e. purely transcendental) over $k$. It is not hard to produce
many examples for which the answer is affirmative. However, no negative
examples were found for more than fifty years.

        The first example of a group for which the Noether Problem has
negative answer is due to Swan \cite{Swan1969}, who proved that over $k=\Q$
the field of invariants of the regular representation of $\Z/47\Z$ is not
rational.

        Over an algebraically closed field $k$ of characteristic zero, the
first negative example is due to Saltman \cite{Saltman1984}; his methods
were subsequently refined by Bogomolov \cite{Bogomolov1988}. The examples
of Saltman and Bogomolov are certain $p$-groups of order $p^9$ and of order
$p^6$, respectively. More recently Hoshi, Kang and Kunyavski\u{\i}
\cite{HKM-noether} found examples of order $p^5$, where $p$ is odd.

        The strategy of Saltman is as follows. If $K/k$ is a field
extension, one can consider the unramified Brauer group $\on{Br}_{nr}(K):=
{H}^2_{nr}(K, \mathbb{Q}/\mathbb{Z}(1))$. It is an abelian group, and it is
trivial if $K/k$ is rational. To find a negative example to the Noether
Problem, Saltman exhibited $G$ and $V$ for which he was able to show that
$\on{Br}_{nr}(k(V)^G)\neq 0$. Later, Bogomolov \cite{Bogomolov1988}
described $\on{Br}_{nr}(k(V)^G)$ purely in terms of the group cohomology of
$G$, as a subgroup of the Schur multiplier $H^3(G,\Z)$. For this reason the
group $\on{Br}_{nr}(k(V)^G)$ is also known as the Bogomolov multiplier of
$G$, and is sometimes denoted by $B_0(G)$.

        \subsection*{Triviality of the motivic class of $BG$}
        In 2009, Ekedahl \cite{ekedahl2009geometric} considered the motivic
class $\set{BG}$ in the Grothendieck ring of stacks $K_0(\on{Stacks}_k)$.
Furthermore, when $k$ has characteristic zero, he used the refined Euler
characteristic introduced in \cite{ekedahl2009geometric} to construct
geometric invariants of $G$, one for each integer $n\geq 1$. For $n=2$, his
construction recovers $B_0(G)$. These invariants are called the Ekedahl
invariants of $G$ in \cite{Martino-TEIFFG, martino2017introduction}, to
which we refer for an overview of this topic.

        We say that $BG$ has trivial motivic class if $\set{BG}=1$ in
$K_0(\on{Stacks}_k)$. Triviality of the class $\{BG\}$ entails triviality
of the Ekedahl invariants and, in particular, of the Bogomolov multiplier;
see \cite[Theorem 5.1]{ekedahl2009geometric}. This implies that in the
aforementioned negative examples to the Noether Problem over an
algebraically closed field of characteristic zero, $\{BG\}$ is not trivial
in $K_0(\on{Stacks}_k)$. The connection between the triviality of
$\set{BG}$ and the Noether Problem is intriguing, but so far remains
largely unexplained.

        A further point of interest in the triviality of $\set{BG}$ comes
from recent work of Totaro \cite{totaro2016motive}, which suggests that it
might be related to five other interesting properties of finite groups:
stable rationality of $BG$, triviality for the birational motive of the
quotient varieties $V/G$, the weak Chow K\"unneth property of $BG$, the
Chow K\"{u}nneth property of $BG$ and the mixed Tate property of $BG$. As
explained by Totaro, it is entirely possible that all these properties are
equivalent, when $k$ is algebraically closed of characteristic zero. For
example, if the Bogomolov multiplier $B_0(G)$ is non-trivial, then all the
above properties fail for $G$.

        We note that over non-algebraically closed fields the above
properties are in general not equivalent. For example, if $k$ is a field of
characteristic zero admitting a biquadratic field extension, there exist
non-constant finite group schemes $G$ over $k$ such that $\set{BG}\neq 1$
and $BG$ is stably rational; see \cite{scavia2018motivic}.

        The known instances of triviality for $\{BG\}$ mainly come from
work of Ekedahl. They are:

        \begin{enumerate}
                \setlength\itemsep{0em}
                \item[--] the group schemes $\mu_n$ of $n$-th roots of
unity, for every $n\geq 1$ (see \cite[Proposition
3.2]{ekedahl2009geometric});
                \item[--] the symmetric groups $S_n$, $n\geq 1$ (see
\cite[Theorem 4.3]{ekedahl2009geometric});
                \item[--]
                all finite subgroups of the group of affine transformations
of $\mathbb{A}^1_{k}$, assuming $k$ algebraically closed (see \cite[p. 8,
Example \textit{ii)} on page 8]{ekedahl2009geometric}).
        \end{enumerate}

        \noindent
        Subsequently, when $k$ is algebraically closed of characteristic
zero, the first author proved triviality of $\{BG\}$ in
$K_0(\on{Stacks}_{k})$ for
        \begin{enumerate}
                \item[--] the finite subgroups of $\on{GL}_3$, \cite[Theorem
2.4]{Martino-TEIFFG}.
        \end{enumerate}
        The proofs of all of these results run along similar lines: one
considers a faithful representation $V$ of $G$, stratifies $V$ according to
the stabilizer, computes $\set{U/G}$ as a polynomial in $\L:=\set{\A^1_k}$,
where $U\c V$ is the open subset where $G$ acts freely, and inductively
computes the other strata.

        \subsection*{Semi-direct products} We now come to the new results
of this paper. We devise a new method for establishing the triviality of
the classes of the classifying stacks of certain finite groups, based on
integral representation theory.

        Let $G$ be a finite constant group, and let $M$ be a $G$-lattice,
that is, $M$ is a finitely generated free abelian group on which $G$ acts
additively. Consider a short exact sequence
\begin{equation}\label{coflasque}
        0\to N\to P\to M\to 0
        \end{equation} where $P$ is a permutation $G$-lattice and $N$ is a
coflasque $G$-lattice. Recall that a $G$-lattice is \emph{permutation} if
it admits a basis which is stable under the $G$-action, and that it is
\emph{coflasque} if $H^1(H,N)=0$ for every subgroup $H$ of $G$. We refer
the reader to \cite[Proposition 1.3]{colliot1987principal} for the
construction of a sequence (\ref{coflasque}).

        We denote by $e(M)$ the \emph{period} of the class of
(\ref{coflasque}) in $\on{Ext}^1_G(M,N)$. The invariant $e(M)$ has a very
interesting geometric interpretation: if $L/F$ is a finite Galois extension
with Galois group $G$, and $U$ is an $F$-torus split by $L$ and whose
character $G$-lattice is isomorphic to $M$, by a theorem of Merkurjev
\cite{merkurjev2010periods} the number $e(M)$ equals the period of a
generic $U$-torsor $X\to\Spec K$, as an element of the group $H^1(K,U)$;
see \cite{merkurjev2010periods} for the precise definitions. In particular,
$e(M)$ does not depend on the choice of (\ref{coflasque}). It is clear that
$e(M)$ is a divisor of $|G|$.

        \begin{thm}\label{mainthm}
                Let $k$ be a field, let $G$ be a finite group, let $T\cong
\G_{\on{m}}^n$ be a split $k$-torus on which $G$ acts multiplicatively, and
let $\hat{T}$ be the character lattice of $T$, viewed as a $G$-lattice via
the induced action. Let $m\equiv\pm 1\pmod{e(\hat{T})}$ be a non-negative
integer, and denote by $T[m]$ the $m$-torsion subgroup of $T$.

                Then we have $\set{B(T[m]\rtimes G)}=\set{BG}$ in
$K_0(\on{Stacks}_{k})$.
        \end{thm}
        In particular, if $\set{BG}=1$ in $K_0(\on{Stacks}_k)$, then
$\set{B(T[m]\rtimes G)}=1$. We note that \Cref{mainthm} makes no
assumptions on the base field $k$: the group scheme $T[m]$ is allowed to be
non-constant and even non-reduced.

        As we have already mentioned, all previous results on the
triviality of $\set{BG}$ have been proved by choosing a faithful
representation $V$ of $G$ and by computing the class of $[V/G]$ using a
suitable stratification. \Cref{mainthm} is of a different nature: its proof
is arithmetic, and makes use of multiplicative invariant theory, the theory
of (non-split) algebraic tori, and Galois cohomology.

        \subsection*{Finite reflection groups }As an application of
\Cref{mainthm}, we prove the triviality of a large number of motivic
classes of classifying stacks of finite complex reflection groups.

        Recall that a finite constant group $G$ is called a \emph{complex
reflection group} if there exists a faithful complex representation $V$ of
$G$ such that $G$ is generated by pseudoreflections (i.e. elements which
fix some complex hyperplane of $V$ pointwise); a useful reference on the
subject is \cite{broue2010introduction}. If $G$ is a complex reflection
group, and $V$ is a reflection representation for $V$, it follows from the
Chevalley-Shephard-Todd Theorem that the Noether Problem for $G$ and $V$
has affirmative answer, and so $BG$ is stably rational. It is then natural
to wonder if some of the other properties listed by Totaro are also true.
In particular, do we have $\set{BG}=1$? If $G=S_n$, then the answer is
affirmative by \cite[Theorem 4.3]{ekedahl2009geometric}, however the proof
does not seem to generalize to other reflection groups.

        For every field $k$ and all integers $m,p,n\geq 1$ such that $p$
divides $m$, consider the $k$-group scheme \[
        G(m,p,n):=\set{(\zeta_1,\dots,\zeta_n,\sigma)\in \mu_m^n\rtimes S_n:
\prod_{i=1}^n\zeta_i^{m/p}=1}.\]
        If $k=\mathbb C$, the $G(m,p,n)$ form the infinite family of
irreducible finite complex reflection groups.
        \begin{cor}\label{gmpn}
                Let $k$ be a field. Assume that $p=1$, or that $p=m\equiv
\pm 1\pmod n$. Then we have $\set{BG(m,p,n)}=1$ in $K_0(\on{Stacks}_k)$.
        \end{cor}
        The proof of \Cref{gmpn} in the case $p=1$ is particularly simple,
and it implies that $\set{BG}=1$ for every finite reflection group $G$ of
type $B_n$; this case may also be deduced from the $Symm$ formalism of
\cite{ekedahl2009geometric}.

        \section*{Acknowledgements} We thank Emanuele Delucchi and Emanuele
Ventura for helpful comments, and the anonymous referee for suggesting
various improvements on the exposition.
        We are grateful to Emanuele Delucchi for his input to earlier
stages of this research project.

        \section{Motivic classes of classifying stacks}
        We begin by recalling that the Grothendieck ring of algebraic
varieties $K_0(\on{Var}_k)$ is the group generated by the isomorphism
classes $\{X\}$ of $k$-schemes of finite type $X$, subject to the relation
$\{X\}=\{Y\}+\{X\setminus Y\}$ for every closed embedding $Y\hookrightarrow
X$. We define a product on $K_0(\on{Var}_k)$ by setting $\{X\}\cdot
\{Y\}:=\{X\times Y\}$ and extending by bilinearity. This makes
$K_0(\on{Var}_k)$ into a commutative ring with identity $1=\set{\Spec k}$.
We denote by $\L$ the class of the affine line $\A^1_k$ in
$K_0(\on{Var}_k)$.

        In \cite{ekedahl2009grothendieck}, Ekedahl constructed a
Grothendieck ring of algebraic stacks as follows.

        \begin{defin}\label{defink0}
                Let $\mc{S}$ be an algebraic stack of finite type over $k$.
The Grothendieck ring of algebraic stacks over $\mc{S}$, denoted
$K_0(\on{Stacks}_{\mc{S}})$, is the quotient of the free abelian group
generated by equivalence classes $\set{\mc{X}}$ of algebraic stacks
$\mc{X}$ finitely presented over $\mc{S}$ and with affine stabilizers by
the following relations:
                \begin{enumerate}
                        \item
$\set{\mc{X}}=\set{\mc{Z}}+\set{\mc{X}\setminus\mc{Z}}$ for every closed
embedding $\mc{Z}\hookrightarrow\mc{X}$;
                        \item
$\set{\mc{E}}=\set{\A^n_{\mc{S}}\times_{\mc{S}}\mc{X}}$ for every vector
bundle $\mc{E}$ of constant rank $n$ over $\mc{X}$.
                \end{enumerate}
                The product on $K_0(\on{Stacks}_{\mc{S}})$ is defined by
$\set{\mc{X}}\set{\mc{Y}}:=\set{\mc{X}\times_{\mc{S}}\mc{Y}}$. The class
$\set{\A^1_{\mc{S}}}$ in $K_0(\on{Stacks}_{\mc{S}})$ is denoted by $\L$, or
by $\L_{\mc{S}}$ if reference to $\mc{S}$ is necessary.
        \end{defin}
        We will be especially interested in the cases, when $\mc{S}=\Spec
k$ or $\mc{S}=BG$ for some finite group scheme $G$ over $k$. There is a
natural ring homomorphism $K_0(\on{Var}_k)\to K_0(\on{Stacks}_k)$, which by
\cite[Theorem 4.1]{ekedahl2009grothendieck} induces an isomorphism
        \[
        K_0(\on{Stacks}_k)\cong K_0(\on{Var}_k)[\mathbb{L}^{-1},
\set{(\mathbb{L}^{n}-1)^{-1}, n\geq 1}].
        \]

        If $f:\mc{S}'\to \mc{S}$ is a morphism of stacks, we get a natural
ring homomorphism $f^*:K_0(\on{Stacks}_{\mc{S}})\to
K_0(\on{Stacks}_{\mc{S}'})$ by pulling back stacks along $f$. In
particular, $f^*(\L_{\mc{S}})=\L_{\mc{S}'}$. If $f$ is finitely presented,
we also have a group homomorphism $f_*:K_0(\on{Stacks}_{\mc{S}'})\to
K_0(\on{Stacks}_{\mc{S}})$, which sends the class of a stack $\mc{X}\to
\mc{S}'$ to the class of $\mc{X}\to \mc{S}'\xrightarrow{f}\mc{S}$. Note
that $f_*$ is not a ring homomorphism, however there is a projection
formula
        \[f_*(f^*C\cdot C')=C\cdot f_*C',\qquad C\in
K_0(\on{Stacks}_{\mc{S}}), \, C'\in K_0(\on{Stacks}_{\mc{S}'}).\] To prove
it, one may reduce to the case when $C=\set{\mc{X}}$ and
$C'=\set{\mc{X}'}$, where $\mc{X}$ is a stack over $\mc{S}$ and $\mc{X}'$
is a stack over $\mc{S}'$, in which case the claim follows from the
$\mc{S}$-isomorphism $\mc{X}_{\mc{S}'}\times_{\mc{S}'}\mc{X}'\cong
\mc{X}\times_{\mc{S}}\mc{X}'$.

        Let $\mc{S}$ be an algebraic stack and let $G$ be a linear
algebraic group scheme over $\mc{S}$, that is, $G$ is flat over $\mc{S}$
and is a closed subgroup of $\on{GL}_{n,\mc{S}}$ for some $n\geq 1$. We
denote by $B_{\mc{S}}G$ the classifying stack of $G$ over $\mc{S}$. If $T$
is a scheme over $\mc{S}$, then by definition $B_{\mc{S}}G(T)$ is the
groupoid whose objects are $G_T$-torsors $P\to T$, and whose arrows are
isomorphisms of $G_T$-torsors over $T$. If $\mc{S}=S$ is a scheme, this is
the usual quotient stack $[S/G]$, where $G$ acts trivially on $S$.
        If $S = \Spec k$, then we relax our notation and simply write $BG$
instead of $B_{\Spec k}G$.

        The structure morphism $B_{\mc{S}}G\to \mc{S}$ sends a $G_T$-torsor
$P\to T$ to its base $T$. The natural morphism $\mc{S}\to B_{\mc{S}}G$,
sending a scheme $T\to \mc{S}$ to the split $G_T$-torsor, is a $G$-torsor.

        The next lemma is an immediate generalization to an arbitrary base
stack $\mc{S}$ of results that are well-known when $\mc{S}=\Spec k$.
        \begin{lem}\label{gloverstacks}
                Let $\mc{S}$ be an algebraic stack of finite type over $k$.
Then:
                \begin{enumerate}[label=(\alph*)]
                        \item\label{gloverstacks1}
$\set{\on{GL}_{n,\mc{S}}}=\prod_{i=0}^{n-1}(\L^n-\L^i)$ in
$K_0(\on{Stacks}_{\mc{S}})$;
                        \item\label{gloverstacks2} if $\mc{X}\to \mc{Y}$ is
a $\on{GL}_{n,\mc{S}}$-torsor of algebraic stacks over $\mc{S}$, then
$\set{\mc{X}}=\set{\on{GL}_{n,\mc{S}}}\set{\mc{Y}}$ in
$K_0(\on{Stacks}_{\mc{S}})$;
                        \item\label{gloverstacks3} we have
$\set{B_{\mc{S}}\on{GL}_{n,\mc{S}}}\set{\on{GL}_{n,\mc{S}}}=1$ in
$K_0(\on{Stacks}_{\mc{S}})$. In particular $\set{\on{GL}_{n,\mc{S}}}$ is
invertible in $K_0(\on{Stacks}_{\mc{S}})$.
                \end{enumerate}
                %If $\mc{S}=[S/\on{GL}_n]$, where $S$ is a scheme, since the
natural map $S\to \mc{S}$ is a $\on{GL}_n$-torsor, we have
$\set{S}=\set{\on{GL}_n}$ in $K_0(\on{Stacks}_{\mc{S}})$ and
$1=\set{\on{GL}_n}\set{{S}}$ in $K_0(\on{Stacks}_{\mc{S}})$. For example,
$\set{\on{Spec} k}=\prod_{i=1}^n(\L^n-\L^i)$ in
$K_0(\on{Stacks}_{B\on{GL}_n})$.
        \end{lem}

        \begin{proof}
                \ref{gloverstacks1} By \cite[Proposition
1.1(i)]{ekedahl2009grothendieck}, the relation
$\set{\on{GL}_{n,k}}=\prod_{i=0}^{n-1}(\L^n-\L^i)$ holds in
$K_0(\on{Stacks}_k)$. The desired formula follows by pulling back along the
structure morphism $\mc{S}\to\Spec k$ of $\mc{S}$.

                \ref{gloverstacks2} Assume first that
$\set{\mc{X}}=\set{\on{GL}_{n,\mc{Y}}}\set{\mc{Y}}$ in
$K_0(\on{Stacks}_{\mc{Y}})$, and denote by $f:\mc{Y}\to \mc{S}$ the
structure morphism. Then, using the projection formula, we deduce

\[f_*\set{\mc{X}}=f_*(\set{\on{GL}_{n,\mc{Y}}}\set{\mc{Y}})=f_*(f^*\set{\on{GL}_{n,\mc{S}}}\set{\mc{Y}})=\set{\on{GL}_{n,\mc{S}}}f_*\set{\mc{Y}},\]
showing that $\set{\mc{X}}=\set{\on{GL}_{n,\mc{S}}}\set{\mc{Y}}$ in
$K_0(\on{Stacks}_{\mc{S}})$. Therefore, we may assume that $\mc{Y}=\mc{S}$,
and the claim becomes $\set{\mc{X}}=\set{\on{GL}_{n,\mc{S}}}$ in
$K_0(\on{Stacks}_{\mc{S}})$. In this case, even though in \cite[Proposition
2.2]{bergh2015motivic} and \cite[Proposition
1.1(ii)]{ekedahl2009grothendieck} it is claimed that
$\set{\mc{X}}=\set{\on{GL}_{n,\mc{S}}}$ in $K_0(\on{Stacks}_k)$, both
proofs show that the equality holds in $K_0(\on{Stacks}_{\mc{S}})$ as well.

                \ref{gloverstacks3} The natural map $\mc{S}\to
B_{\mc{S}}\on{GL}_{n,\mc{S}}$ is a $\on{GL}_{n,\mc{S}}$-torsor, hence the
conclusion follows from \ref{gloverstacks2}.
        \end{proof}

        The next lemma will be essential to derive the results of
\Cref{groupsofmulttype}. It will allow us to reduce questions about classes
in $K_0(\on{Stacks}_{BG})$ of representable stacks over $BG$ to questions
about schemes. %; see \Cref{groupsofmulttype}.

        \begin{lem}\label{enoughonspaces}
                Let $k$ be a field, and let $G$ be a linear algebraic group
over $k$. Then there exist a $k$-variety $X$ and a morphism $f:X\to BG$ such
that the pullback map $f^*:K_0(\on{Stacks}_{BG})\to K_0(\on{Stacks}_X)$ is
injective.
        \end{lem}

        \begin{proof}
                Choose an embedding of $G\hookrightarrow \on{GL}_n$ for
some $n\geq 1$, and set $X:=\on{GL}_n/G$. The variety $X$ is a homogeneous
space under $\on{GL}_n$ and has a $k$-point with stabilizer isomorphic to
$G$, hence $[X/\on{GL}_n]\cong BG$. The canonical map $f:X\to
[X/\on{GL}_n]=BG$ is a $\on{GL}_n$-torsor. For a stack $\phi:\mc{X}\to BG$,
let $f^*\mc{X}:=\mc{X}\times_{\phi,BG,f}X$, so that
$f^*\set{\mc{X}}=\set{f^*\mc{X}}$ in $K_0(\on{Stacks}_X)$. The first
projection $f^*\mc{X}\to\mc{X}$ is a $\on{GL}_n$-torsor, hence by
\Cref{gloverstacks}\ref{gloverstacks2}
\begin{equation}\label{pushpull}f_*f^*\set{\mc{X}}=f_*\set{f^*\mc{X}}=\set{\mc{X}}\set{\on{GL}_n}\end{equation}
in $K_0(\on{Stacks}_{BG})$.

                Let $C\in K_0(\on{Stacks}_{BG})$. We may write
$C=\sum_i\set{\mc{X}_i}-\sum_j\set{\mc{Y}_j}$, for some algebraic stacks
$\mc{X}_i$ and $\mc{Y}_j$ over $BG$. Using (\ref{pushpull}) on each term,
we get
\[f_*f^*C=\sum\set{\mc{X}_i}\set{\on{GL}_n}-\sum\set{\mc{Y}_i}\set{\on{GL}_n}=\sum(\set{\mc{X}_i}-\set{\mc{Y}_i})\set{\on{GL}_n}=C\set{\on{GL}_n}\]
in $K_0(\on{Stacks}_{BG})$. In other words,
$f_*f^*:K_0(\on{Stacks}_{BG})\to K_0(\on{Stacks}_{BG})$ is the map of
multiplication by $\set{\on{GL}_n}$. By
\Cref{gloverstacks}\ref{gloverstacks3}, $\set{\on{GL}_n}$ is invertible,
hence $f^*$ is injective.
        \end{proof}

        %\section{Preliminaries on (\ref{expected})}
        \medskip
        Let $S$ be a scheme, and let $G$ be a linear algebraic group scheme
over $S$. We say that $G$ is \emph{special} if
$H^1_{\operatorname{fppf}}(T,G)=H^1_{\operatorname{Zar}}(T,G)$ for every
$S$-scheme $T$, that is, if $G$-torsors are Zariski-locally trivial over
any $S$-scheme. The following result has been proved in
\cite{bergh2015motivic} in the case, when $S=\Spec k$ for some field $k$.

        \begin{lem}\label{special}
                Let $S$ be a scheme and $G$ be a special $S$-group.
                \begin{enumerate}[label=(\alph*)]
                        \item\label{special1} If $\pi:X\to Y$ is a
$G$-torsor of $S$-schemes, then $\set{X}=\set{G}\set{Y}$ in
$K_0(\on{Stacks}_S)$.
                        \item\label{special2} If $\mc{X}\to\mc{Y}$ is a
$G$-torsor of $S$-stacks, then $\set{\mc{X}}=\set{G}\set{\mc{Y}}$ in
$K_0(\on{Stacks}_S)$.
                        \item\label{special3} We have $\set{B_SG}\set{G}=1$
in $K_0(\on{Stacks}_S)$.
                \end{enumerate}
        \end{lem}

        \begin{proof}
                \ref{special1} The proof of \cite[Proposition
2.1]{bergh2015motivic} immediately generalizes; we include it for
completeness. Since $G$ is special, there exists a non-empty open subscheme
$X_1$ of $X$ such that $Y_1:=\pi^{-1}(X_1)$ is isomorphic to
$X_1\times_SG$. Iterating this procedure, we eventually obtain a
stratification $X=\amalg X_i$, such that $Y_i:=\pi^{-1}(X_i)\cong
X_i\times_SG$. By the scissor relation, we have $\set{X}=\sum\set{X_i}$ and
$\set{Y}=\sum \set{Y_i}$. Since $\set{Y_i}=\set{X_i}\set{G}$ for every $i$,
we conclude that \[\set{Y}=\sum \set{Y_i}= \sum\set{X_i}\set{G}=
                \left(\sum\set{X_i}\right)\set{G}=\set{X}\set{G}.\]

                \ref{special2} %Using \Cref{gloverstacks} instead of
\cite[Proposition 2.2]{bergh2015motivic}, the proof of \cite[Proposition
2.3]{bergh2015motivic} adapts without difficulties.
                Assume that $\mc{Z}$ is a closed substack of $\mc{Y}$, with
open complement $\mc{U}$, and that the claim holds for
$\mc{X}_{\mc{Z}}\to\mc{Z}$ and $\mc{X}_{\mc{U}}\to\mc{U}$. Then, by the
scissor relation

\[\set{\mc{X}}=\set{\mc{X}_{\mc{U}}}+\set{\mc{X}_{\mc{Z}}}=\set{\mc{U}}\set{G}+\set{\mc{Z}}\set{G}=\set{\mc{Y}}\set{G}.\]
By noetherian induction, it is enough to show the claim for a non-empty
open substack of $\mc{X}$. By \cite[Proposition 3.5.9]{kresch1999cycle},
$\mc{Y}$ is stratified by stacks of the form $[U/\on{GL}_{n,S}]$, where $U$
is a scheme over $S$. Hence we may assume that
$\mc{Y}=[Y/\on{GL}_{n,{S}}]$, where $n\geq 1$ and $Y$ is a scheme. We have
a cartesian diagram
                \begin{equation*}
                \begin{tikzcd}
                \mc{X}_{Y} \arrow[d] \arrow[r,"\phi"] & \mc{X} \arrow[d]\\
                Y \arrow[r,"f"] & \mc{Y}
                \end{tikzcd}
                \end{equation*}
                where the horizontal maps are $\on{GL}_{n,{S}}$-torsors. It
follows from \Cref{gloverstacks}\ref{gloverstacks2} that
$\set{Y}=\set{\on{GL}_{n,{S}}}\set{\mc{Y}}$ and
$\set{\mc{X}_Y}=\set{\on{GL}_{n,{S}}}\set{\mc{X}}$. Since $Y$ is a scheme,
by \ref{special1} we have $\set{\mc{X}_Y}=\set{\on{GL}_{n,{S}}}\set{Y}$.
Combining these relations, we arrive to
$\set{\on{GL}_{n,{S}}}(\set{\mc{X}}-\set{G}\set{\mc{Y}})=0$. By
\Cref{gloverstacks}\ref{gloverstacks3}, $\set{\on{GL}_{n,{S}}}$ is
invertible, hence $\set{\mc{X}}-\set{G}\set{\mc{Y}}=0$.

                \ref{special3} If $T$ is a scheme and $T\to B_SG$ is a
morphism, corresponding to a $G$-torsor $P\to T$, by \ref{special1} we have
$\set{P}=\set{G}\set{T}$. Applying \ref{special2} to $S\to B_SG$ and
$C=\set{G}$, we get $1=\set{S}=\set{B_SG}\set{G}$ in $K_0(\on{Stacks}_S)$.
        \end{proof}

        \section{Groups of multiplicative type over
stacks}\label{groupsofmulttype}
        Let $G$ be a finite group and $M$ be a $G$-module of rank $n$. We
define $\mc{T}_M:=[\Spec k[M]/G]$ and $\mc{B}_M:=B(\Spec k[M]\rtimes G)$.

        The representable morphism $\mc{T}_M\to BG$ exhibits $\mc{T}_M$ as
a group object over $BG$. If $f:X\to BG$ is a morphism, corresponding to a
$G$-torsor $\pi:Y\to X$, then $T_M:=X\times_{f,BG,\pi}\mc{T}_M$ is a group
of multiplicative type split by $Y$. It is obtained by twisting $\Spec
k[M]\times X$ by $\pi$ (using the $G_X$-action on $\Spec k[M]\times X$). We
will be mostly interested in the case, when $M$ is torsion-free. In that
case $T_M$ is a relative torus over $X$.

        If $\cl{x}$ is a geometric point of $X$, then $\pi$ corresponds to
a surjection $\pi_1(X,\cl{x})\to G$, and $T_M$ corresponds to the
representation of $\pi_1(X,\cl{x})\to G\to \on{GL}(M)$.

        The projection $\mc{B}_M\to BG$ admits a section, making $\mc{B}_M$
into a neutral gerbe over $BG$. Its fibers are the twists of $\Spec k[M]$
using the $G$-action, as we now explain.

        More generally, let us consider the semidirect product $N\rtimes H$
of two linear algebraic $k$-groups $N$ and $H$. If $f:X\to BH$ is a
morphism, corresponding to an $H$-torsor $\pi:Y\to X$, we have a cartesian
diagram
        \begin{equation}\label{semidirect}
        \begin{tikzcd}
        \text{$[X/\prescript{\pi}{}{N}]$} \arrow[d] \arrow[r,"\phi"] &
B(N\rtimes H) \arrow[d]\\
        X \arrow[r,"f"] & BH.
        \end{tikzcd}
        \end{equation}
        Here $\prescript{\pi}{}{N_X}$ is the twist of $N_X$ by the torsor
$\pi:Y\to X$, using the $H_X$-action. It is a group scheme over $X$. The
morphism $\phi$ is constructed in the following way: twist the inclusion
$N_X\to (N\rtimes H)_X$ by the $(N\rtimes H)_X$-torsor $Y\times^H(N\rtimes
H)\to X$ to obtain a morphism $\prescript{\pi}{}{N_X}\to
\prescript{\pi}{}{(N\rtimes H)_X}$. Since $\prescript{\pi}{}{(N\rtimes
H)_X}$ is an inner form of $(N\rtimes H)_X$, we have
$[X/\prescript{\pi}{}{(N\rtimes H)_X}]\cong [X/(N\rtimes H)_X]$. We define
$\phi$ as the composition \[[X/\prescript{\pi}{}{N}]\to
[X/\prescript{\pi}{}{(N\rtimes H)_X}]\cong [X/(N\rtimes H)_X]\cong X\times
B(N\rtimes H)\xrightarrow{\on{pr}_2} B(N\rtimes H).\]
        Specialization of diagram (\ref{semidirect}) to our situation
yields a cartesian diagram
        \begin{equation}\label{semidirect2}
        \begin{tikzcd}
        B_XT_M \arrow[d] \arrow[r] & \mc{B}_M \arrow[d]\\
        X \arrow[r,"f"] & BG
        \end{tikzcd}
        \end{equation}
        for every morphism $f:X\to BG$.

        Recall that, by definition, an invertible $G$-lattice $M$ is a
direct summand of a permutation lattice \cite{MR1027091,
lorenz2006multiplicative}.
        \begin{lem}\label{permutation}
                Let $M$ be an invertible $G$-lattice (for example, a
permutation $G$-lattice). Then:
                \begin{enumerate}[label=(\alph*)]
                        \item\label{permutation0'} for every morphism
$f:X\to BG$, the $X$-torus $f^*\mc{T}_M$ is special;
                        \item\label{permutation0}
$\set{\mc{T}_M}\set{\mc{B}_M}=1$ in $K_0(\on{Stacks}_{BG})$;
                \end{enumerate}
        \end{lem}

        \begin{proof}
                Let $P$ be a permutation $G$-lattice such that $M$ is a
direct summand of $P$. Let $Y\to X$ be the $G$-torsor corresponding to $f$.
Then $f^*\mc{T}_P=R_{Y/X}(\G_{\on{m},Y})$, $f^*\mc{B}_P=
B_X(R_{Y/X}(\G_{\on{m},Y}))$ and $f^*\mc{T}_M$ is a direct factor of
$f^*\mc{T}_P$.

                \ref{permutation0'} By Shapiro's lemma for \'etale
cohomology \cite[Lemma 29.6]{knus1998involutions}, the $X$-torus
$R_{Y/X}(\G_{\on{m},Y})$ is special. Since $f^*\mc{T}_M$ is a direct factor
of $f^*\mc{T}_P$, the cohomology of $f^*\mc{T}_M$ is a direct summand of
the cohomology of $f^*\mc{T}_P$, hence $f^*\mc{T}_M$ is also special.

                For the rest of the proof we fix a $k$-variety $X$ and a
morphism $f:X\to BG$ such that $f^*:K_0(\on{Stacks}_{BG})\to
K_0(\on{Stacks}_X)$ is injective, by \Cref{enoughonspaces}.

                \ref{permutation0} By \Cref{special}, we obtain that
\[f^*(\set{\mc{T}_M}\set{\mc{B}_M})=\set{f^*\mc{T}_M}\set{f^*\mc{B}_M}=\set{R_{Y/X}(\G_{\on{m},Y})}\set{B(R_{Y/X}(\G_{\on{m},Y}))}=1\]
in $K_0(\on{Stacks}_X)$. Since $f^*$ is injective, we conclude that
$\set{\mc{T}_M}\set{\mc{B}_M}=1$.
        \end{proof}

        \begin{prop}\label{stackytori}
                Let \[0\to M'\to M\to M''\to 0\] be a short exact sequence
of $G$-modules. If either $M$ or $M''$ is an invertible $G$-lattice, then
\[\set{\mc{T}_{M'}}=\set{\mc{T}_M}\set{\mc{B}_{M''}}\] in
$K_0(\on{Stacks}_{BG})$.
        \end{prop}

        \begin{proof}
                \Cref{enoughonspaces} gives a $k$-variety $X$ and a
morphism $f:X\to BG$ such that $f^*:K_0(\on{Stacks}_{BG})\to
K_0(\on{Stacks}_X)$ is injective. Fix a geometric point $\cl{x}$ of $X$.
The morphism $f$ corresponds to a $G$-torsor $Y\to X$, and as discussed
this in turn corresponds to a surjection $\pi_1(X,\cl{x})\to G$. The
sequence $0\to M'\to M\to M''\to 0$ is then a sequence of integral
$\pi_1(X,\cl{x})$-representations, hence we obtain a short exact sequence
of algebraic tori over $X$:
                \[1\to T_{M''}\to T_M\to T_{M'}\to 1.\]
                Since $f^*$ is injective, it is enough to show that
$f^*\set{\mc{T}_{M'}}=f^*(\set{\mc{T}_M}\set{\mc{B}_{M''}})$, that is
                \begin{equation}\label{stackytorisuffices}
                \set{T_{M'}}=\set{T_M}\set{B_XT_{M''}} \text{ in
$K_0(\on{Stacks}_X)$.}
                \end{equation}
                If $M''$ is invertible, by
\Cref{permutation}\ref{permutation0} $T_{M''}$ is special. Since $T_M\to
T_{M'}$ is a $T_{M''}$-torsor, (\ref{stackytorisuffices}) follows from
\Cref{permutation}.

                If $M$ is invertible, then (\ref{stackytorisuffices})
follows from \cite[Proposition 2.9]{bergh2015motivic}; we repeat the
argument here. The map $T_{M'}\to [{T}_{M'}/T_M]$ is a $T_{M}$-torsor.
Since $T_M$ acts transitively on $T_{M'}$ and $T_{M'}$ admits an $X$-point
with stabilizer $T_{M''}$ (for example, the identity section), we have
$[T_{M'}/T_M]\cong B_XT_{M''}$. By \Cref{permutation}\ref{permutation0},
$T_M$ is special. It follows that
$\set{T_{M'}}=\set{T_M}\set{[{T}_{M'}/T_M]}=\set{T_M}\set{B_XT_{M''}}$.
        \end{proof}

        \section{Proof of Theorem \ref{mainthm}}

        In this section, we prove the main results of this paper. The proof
of \Cref{mainthm} was inspired by the second author's previous work
\cite{scavia2018noether}, where it was shown that if $T$ is a non-split
$k$-torus, the stable rationality of $BT[n]$ has an answer that is periodic
in $n$, with period dividing the period of the generic $T$-torsor; see
\cite{merkurjev2010periods} for the definition. The period of the generic
$T$-torsor equals the number $e(M)$ that was defined in the Introduction.

        \begin{proof}[Proof of \Cref{mainthm}]
                Let $M$ be the character $G$-lattice of $T$. Consider the
following diagram with exact rows and columns
                \begin{equation}\label{bigdiagram}
                \begin{tikzcd}
                &&& 0 \arrow[d] \\
                & 0 \arrow[d] && M \arrow[d] \\
                0 \arrow[r] & N \arrow[r] \arrow[d] & P \arrow[r] \arrow[d]
& M \arrow[r] \arrow[d] & 0 \\
                0 \arrow[r] & N_m \arrow[r] \arrow[d] & P \arrow[r]  & M/mM
\arrow[r] \arrow[d] & 0 \\
                & M \arrow[d] && 0 \\
                & 0
                \end{tikzcd}
                \end{equation}
                where the first row is (\ref{coflasque}). Denote by
$\alpha$ the class of (\ref{coflasque}) in $\on{Ext}^1_G(M,N)$. By
\cite[Proposition 3.1]{scavia2018noether} the class of \[0\to N\to N_m\to
M\to 0\] in $\on{Ext}^1_G(M,N)$ is $m\alpha$. The order of $\alpha$ in
$\on{Ext}^1_G(M,N)$ divides $e(M)$ (it actually equals $e(M)$, by
\cite[Theorem 3.1]{merkurjev2010periods}). It follows that $N_m\cong
N_{m'}$ when $e(M)\mid m-m'$. Recall that if $\gamma\in \on{Ext}^1_G(M,N)$
is the class of \[0\to N\to Q\xrightarrow{\eta} M\to 0,\] then $-\gamma$ is
represented by \[0\to N\to Q\xrightarrow{-\eta} M\to 0.\] It follows that
$N_m\cong N_{m'}$ when $e(M)\mid m+m'$ as well.

                Assume that $m\equiv \pm 1 \pmod{e(M)}$. Then $N_m\cong
N_1\cong P$ and the second row of (\ref{bigdiagram}) becomes \[0\to P\to
P\to M/mM\to 0.\] Since $P$ is a permutation $G$-lattice, by
\Cref{stackytori} we deduce that $\set{B(T[m]\rtimes G)}=1$ in
$K_0(\on{Stacks}_{BG})$. This implies that $\set{B(T[m]\rtimes
G)}=\set{BG}$ in $K_0(\on{Stacks}_{k})$, as desired.
        \end{proof}

        We now prove our main result concerning reflection groups
$G(m,p,n)$. As we anticipated in the introduction, the case $p=1$ contains
all groups of type $B_n$ and is particularly easy.

        \begin{proof}[Proof of \Cref{gmpn}]
                Denote by $U_n$ the $S_n$-lattice $\Z[S_n/S_{n-1}]$, and let
$e_1,\dots,e_n$ be the standard basis of $U_n$.

                Assume first that $p=1$. We have
$BG(m,1,n)=\mc{B}_{U_n/mU_n}$. Applying \Cref{stackytori} to the short
exact sequence \[0\to U_n\xrightarrow{\times m} U_n\to U_n/mU_n\to 0,\] we
obtain $\set{BG(m,1,n)}=1$ in $K_0(\on{Stacks}_{BS_n})$. Since
$\set{BS_n}=1$ in $K_0(\on{Stacks}_k)$ by \cite[Theorem
4.3]{ekedahl2009geometric}, we conclude that $\set{BG(m,1,n)}=1$ in
$K_0(\on{Stacks}_k)$.

                Assume now that $p=m$. Let $A_{n-1}$ be the $S_n$-module
defined by the short exact sequence
                \begin{equation}\label{define}0\to \Z\xrightarrow{\epsilon}
U_n\to A_{n-1}\to 0,\end{equation} where $\epsilon(1)=\sum_{i=1}^n e_i$.
Let $S_n$ act on $\G_{\on{m}}^n$ by permutation of the coordinates, and
consider the subtorus of $\G_{\on{m}}^n$ given by
\[T_n:=\set{(x_1,\dots,x_n)\in \G_{\on{m}}^n: \prod_{i=1}^nx_i=1}.\] The
character $S_n$-lattice of $T_n$ is $A_{n-1}$, and $G(m,m,n)=T_n[m]\rtimes
S_n$. In view of \Cref{mainthm}, to complete the proof it is enough to show
that $e(A_{n-1})= n$. Since (\ref{define}) is a coflasque resolution of
$A_{n-1}$, it is enough to show that the class of (\ref{define}) in
$\on{Ext}_{S_n}^1(A_{n-1},\Z)$ has order $n$. This is equivalent to showing
that the class of the dual of (\ref{define}) has order $n$ in
$\on{Ext}_{S_n}^1(\Z,A_{n-1}')$, where $A_{n-1}'$ is the $S_n$-lattice dual
to $A_{n-1}$. This is well-known; see e.g. \cite[Example
4.1]{merkurjev2010periods}. It follows from \Cref{mainthm} that
$\set{BG(m,m,n)}=\set{BS_n}$ in $K_0(\on{Stacks}_k)$ when $m\equiv \pm
1\pmod n$. We have $\set{BS_n}=1$ by \cite[Theorem
4.3]{ekedahl2009geometric}, and the conclusion follows.
        \end{proof}

        We conclude by stating the following conjecture.

        \begin{conj}
                Let $W$ be a complex reflection group. Then $\{BW\}=1$ in
$K_0(\on{Stacks}_{\mathbb C})$.
        \end{conj}
        This statement was claimed as the main theorem of the first online
version of the preprint \cite{delucchi2015subspace} of Emanuele Delucchi and
the first author, but its proof was later shown to be defective.

\end{document}